\documentclass[11pt]{article}
\usepackage{latexsym}
\usepackage{amsfonts}
\usepackage{epsfig}

\makeatletter

\def\@footnotetext#1{\insert\footins{%
\footnotesize
    \interlinepenalty\interfootnotelinepenalty
    \splittopskip\footnotesep
    \splitmaxdepth \dp\strutbox \floatingpenalty \@MM
    \hsize\columnwidth \@parboxrestore
   \edef\@currentlabel{\csname p@footnote\endcsname\@thefnmark}\@makefntext
    {\rule{\z@}{\footnotesep}\ignorespaces
      #1\strut}}}
\def\abstract{\small\quotation{\hskip-\parindent\sc Abstract.}}

\def\classification{\@ifnextchar [{\@xfootnotenext}%
   {\begingroup\let\protect\noexpand
      \xdef\@thefnmark{}\endgroup
    \@footnotetext}}

\begin{document}
\classification {{\it 2000 Mathematics
Subject Classification:} Primary 20E05, 20F05; Secondary  57M05,
57M20. } 

\title {}

\begin{center}

\end{center} 

\vskip 3.3cm

\begin{center}
{\bf \Large    On the Andrews-Curtis equivalence} 

\bigskip

{\bf Alexei D. Myasnikov \hskip 20pt Alexei G. Myasnikov \\
    Vladimir Shpilrain}

\end{center} 
\bigskip

\begin{abstract}
\noindent  The Andrews-Curtis conjecture claims that every balanced 
presentation of the trivial group can be reduced to the  standard one by a
sequence of ``elementary  transformations" which are Nielsen 
transformations augmented by arbitrary conjugations. It is a  prevalent  
opinion that this conjecture is false;  however, not   many
potential counterexamples are known. In this paper, we show that some 
of the previously proposed examples are actually not counterexamples. 
We hope that the tricks  we used in constructing relevant chains of
elementary transformations will be useful to  those who attempt to 
establish the Andrews-Curtis equivalence in other situations. 

 On the other hand, we give two rather general and simple methods for 
 constructing balanced presentations of the trivial group; 
some of these presentations can be considered potential counterexamples to 
the Andrews-Curtis conjecture.
 One of the methods is based on a simple combinatorial idea of composition 
of group presentations, whereas the other one uses ``exotic" knot diagrams 
of the unknot. 
 
  We also consider the Andrews-Curtis equivalence in  
metabelian groups and reveal some interesting connections 
of relevant problems to well-known problems in $K$-theory.

\end{abstract}

\bigskip

\noindent {\bf 1. Introduction }
\bigskip

  Let  $F = F_n$   be the free group of a finite rank  
$n \ge 2$  with a  set $X = {\{}x_1, ..., x_n {\}}$  of
free generators. Let   $Y  = {\{}y_1, ..., y_m {\}}$  
 be an arbitrary
finite tuple of elements of the  group $F$.  Consider  the  following 
elementary  transformations  that  can  be applied to  $Y$: 
\smallskip 

\noindent {\bf (AC1)} $y_i$   is replaced by  $y_i y_j$ for some 
 $j \ne i$;
\smallskip 

\noindent {\bf (AC2)} $y_i$ is replaced by $y_i^{-1}$;
\smallskip 

\noindent {\bf (AC3)} $y_i$ is replaced by  $g y_i g^{-1}$ for some 
$g \in F$. 
\smallskip 

    It is understood that $y_j$  does not change if  $j \ne i$. 

    We say that two tuples  $Y$  and  $\widetilde Y$ are  Andrews-Curtis 
equivalent  (AC-equivalent) if one of them can be obtained from another
by  applying a sequence of transformations  (AC1)--(AC3). 

 Denote by $\langle Y
\rangle^F$ the normal closure of a set $Y$ in $F$. The following conjecture appears to
be of interest in topology as well as group theory (see \cite{AC}):
\smallskip

\noindent {\bf The Andrews-Curtis conjecture.} Let $Y  = 
{\{}y_1, ..., y_n {\}}$ be a tuple of elements of $F = F_n.$ Then $\langle Y
\rangle^F = F$  if and only if $Y$ is AC-equivalent to
$X= {\{}x_1, ..., x_n {\}}$.
\smallskip 

 This also has an obvious interpretation in  a different 
 language: every balanced presentation of the trivial group can be
reduced to the  standard one by a sequence of transformations 
(AC1)--(AC3). We  shall say that two presentations of a group are 
AC-equivalent if the sets of relators in these presentations are 
AC-equivalent. 

    A weaker form of this conjecture (which is of  greater  interest in 
topology) arises when one allows two more elementary transformations: 
\smallskip 

\noindent {\bf (AC4)} $Y$ is extended to  ${\{}y_1, ..., y_m,
x_{\nu}{\}}$, where $x_{\nu}$ is a new free generator (i.e., $y_1, ..., y_m$ 
do not depend on  $x_{\nu}$). 

\noindent {\bf (AC5)} the converse of (AC4). 
\smallskip 

 If  $\widetilde Y$  can  be  obtained  from $Y$  by  means of
transformations (AC1--AC5), we say that $\widetilde Y$ is {\it stably} 
AC-equivalent to $Y$. 
\smallskip 

\noindent {\bf The weak Andrews-Curtis conjecture.} For a tuple  
$Y  = {\{}y_1, ..., y_n {\}}$, one has $\langle Y\rangle^F  = F$ 
if  and  only  if $Y$ is stably  AC-equivalent to $X= {\{}x_1, ..., x_n {\}}$. 
\smallskip 

    We mention here a connection of  these conjectures  to 
2-dimensional  homotopy  theory.   There  is  a  standard   way   of  
constructing   a  2-dimensional cell   complex from an arbitrary 
finite presentation of a    group (see e.g. \cite[p.117]{LS}).  
(A complex like that has a single vertex; the number of 1-dimensional 
cells (edges) equals the number of generators in a given group 
presentation,  and the number of 2-dimensional cells equals the number of
defining relators).  Two  cell complexes are called simple-homotopy 
equivalent if one of them can be obtained from another by  a sequence of
elementary expansions and collapses (where ``elementary" means that only
one cell  at a time is added/cancelled). 

 It is well-known  (see  e.g. \cite{Wr} for  details) that (stably) 
AC-equivalent sets of relators give rise  to  simple-homotopy  equivalent
complexes, since the  transformations (AC1-AC5) applied to defining
relators of the fundamental group of a 2-complex, are mimicked by
elementary expansions and collapses applied to this 2-complex. The point
however  is that those elementary expansions/collapses that mimic the
AC-transformations, only involve cells of dimension 3 or lower, whereas
to establish a simple-homotopy  equivalence of two 2-complexes, one is
allowed to use cells of arbitrary dimensions. 
 
 Now, it is known that 
any  two 2-complexes with the trivial fundamental group 
(i.e., 2-complexes contractible to a point) are simple-homotopy 
equivalent, and by a result of Whitehead \cite{Wh}, one does not need to use cells
of  dimension higher than 4 to establish a simple-homotopy  equivalence 
in that situation. 

 So, the point of the (topological) Andrews-Curtis conjecture is that one
does not actually have to use cells of dimension higher than 3  when
treating contractible 2-complexes.  This is equivalent to the weak Andrews-Curtis 
conjecture above (see \cite{Wr}). 
\smallskip 

 We also note that the Andrews-Curtis conjecture admits several   
generalizations -- see the surveys \cite{GK}, \cite{Hog1},  \cite{Hog}, 
\cite{Metzler2}. 
 Here we would like to point out the following very natural  
generalization. 
\smallskip 

Let $G$ be a group and $k$  a positive integer. We say that $G$ 
{\it satisfies the Andrews-Curtis property $ACP_k$}
if for any two $k$-tuples $U, V \in G^k$ the following equivalence
holds:

$$\langle U\rangle^G = \langle V\rangle^G = G \Longleftrightarrow U \sim V,$$
where  $U \sim   V$  means AC-equivalence. (Note that  the definition of
 AC-equivalence can be carried on to arbitrary groups.)

\medskip
\noindent {\bf Question}. {\it Is it true that for every group $G$, there
is an integer $k_0$ such that $G$ satisfies $ACP_k$  for every
$k \geq k_0$ ?} 

\smallskip
 Note that if $A$ is a finite abelian group represented as 
$$
A \simeq {\bf Z}_{e_1} \times\cdots\times
{\bf Z}_{e_r},
$$
where $e_1 \mid e_2 \mid \ldots \mid e_r$, then $A$ does not
satisfy $ACP_r$,  but it does satisfy $ACP_{k}$  for every $k > r$. 
This result is due to B. H. Neumann and H. Neumann \cite{nn}  (see also
 \cite{dg}).

 For some groups, the above question is known to have  the affirmative answer 
 (see \cite{Myasnikov} or
 discussion below), but we are unaware of any example of a  
 group where the answer  would be negative. 
 
\smallskip 

 It is a prevalent opinion that the Andrews-Curtis conjecture is false; 
however, not   many potential counterexamples are known, the reason 
being the difficulty of constructing non-trivial (balanced) 
presentations of the trivial group. Here we give a (probably, partial) 
list of previously known potential counterexamples: 
\smallskip 

\noindent {\bf (1)} $\langle x,y| ~x^{-1}y^2 x = y^3, ~y^{-1}x^2 y = x^3\rangle$.

\noindent {\bf (2)} $\langle x,y, z| ~y^{-1}x y = x^2, ~z^{-1}y z = y^2, 
~x^{-1}z x = z^2\rangle$. 

\noindent {\bf (3)} $\langle x,y| ~x^4y^3 = y^2x^2, x^6y^4 = y^3x^3 \rangle$. 

\noindent {\bf (4)} $\langle x,y| ~x^n = y^{n+1}, ~xyx=yxy\rangle,  ~n \ge 2$.

\smallskip 

 The first two examples are probably the most established ones; see 
\cite{Olga} for discussion.   The presentation (3) is attributed to 
M.Wicks (by folklore), but we do not have a relevant reference. The series (4),
 which comes from \cite{AK},
 is  valuable since it provides  an infinite
collection of examples.  The 
following   series has appeared   recently \cite{MS}: 
\smallskip 

\noindent {\bf (5)} $\langle x,y| ~x^{-1}y^n x = y^{n+1}, ~x=w\rangle$, ~where 
$n \ge 1$, and  $w$  is a word in $x$ and  $y$ with exponent sum 0 on $x$. 
\smallskip 

 Note   that the series (5) contains the example (1), and, furthermore,  
 Proposition 1.3  
 below shows that every presentation in the series (4) is AC-equivalent 
 to a presentation in the series (5). 

 Finally, the following series, due to C. Gordon (see \cite{Brown}), was 
brought to our attention by W. Metzler.

\smallskip 

\noindent {\bf (6)} $\langle x,y| ~x=[x^m, y^n], ~y=[x^p, y^q]\rangle$, ~where 
$m, n, p, q$ ~are arbitrary integers. 

\smallskip

 It is natural to start out by trying  potential counterexamples with minimal
total length of relators.   Assisted by a computer, we were able to ``crack" 
all potential counterexamples with  total length of up to 12. More precisely 
(see \cite{MMS}): any presentation of the trivial
group of the form $\langle x,y \mid r(x,y), s(x,y) \rangle$
~with  ~$|r(x,y)| + |s(x,y)| \leq 12$  ~is AC-equivalent to the
standard one. 

 We also mention here an auxiliary result of independent interest (see \cite{MMS}):
if $G$ is a group defined by a presentation
   ~$\langle \ x,y \mid r(x,y),\  s(x,y) \ \rangle$, where
$|r(x,y)| + |s(x,y)| \leq 12$, 
~and  the abelianization of $G$ is trivial, then $G$ is either the trivial group
or   isomorphic to the following finite group of order 120: 
$~\langle x,y \mid y x y = x^2, \ x y x = x^4\  \rangle.$

These results were obtained by using  MAGNUS 
software package for symbolic computation.
\smallskip 

 In this paper, we  treat a couple of most interesting (in our 
 opinion) presentations covered by the result  of \cite{MMS}
  cited above, because we 
hope that the tricks  we used in constructing relevant chains of
elementary transformations will be useful to  those who attempt to 
establish the Andrews-Curtis equivalence in other situations:
\smallskip 

 \noindent {\bf Proposition 1.1.} \cite{MMS} Each of the  following presentations 
of the trivial group is AC-equivalent to $\langle x,y| ~x, y\rangle$: 
\smallskip 

 \noindent {\bf (a)} $\langle x,y| ~x^2 = y^3, ~xyx=yxy\rangle$. 
\smallskip 

\noindent {\bf (b)} $\langle x,y| ~x^{-1}yx =y^2, ~x=yx^2yx^{-2}\rangle$.  

 \noindent {\bf (c)} $\langle x,y| ~x^{-1}y^2 x = y^3, ~x^2=yxy^{-1}\rangle$. 
\smallskip 

\noindent {\bf (d)} $\langle x,y| ~x^{-1}y^2 x = y^3, ~x^2 = yxy\rangle$. 
\smallskip

We  note that some time ago, S. Gersten (unpublished) showed that 
the presentation (a) is {\it stably} AC-equivalent to 
$\langle x,y| ~x, y\rangle$. 
\smallskip 

  Thus, the minimal total length of relators in potential counterexamples 
  that still stand, is 13,  
as in 
$$\langle x,y| ~x^3 = y^4, ~xyx=yxy\rangle.$$

 \indent Recently, Havas and Ramsay  \cite{Havas} showed that  this is, in fact, 
the only (up to AC-equivalence) possible counterexample of length 13.  
 We were not able 
  to crack this example, but the following fact might be of interest: 
\smallskip 
 
 \noindent {\bf Proposition 1.2.} 
 The presentation $\langle x,y| ~x^3 = y^4, ~xyx=yxy\rangle$ 
~is stably AC-equivalent to 
 $\langle x,y | ~x^4 = y x^2 y^{-1} x^{-1} y x^2 y^{-1}, 
~y = [x^2, y]^3 \rangle$. (Our commutator notation is $[a, b] = ab a^{-1} b^{-1}$.)
\smallskip 

This gives what seems to be  the first 
 example of two stably AC-equivalent sets in a free group,  that are not 
 known to be AC-equivalent. The longer presentation   eludes all 
 attempts  to decrease the total length of relators by using the 
 transformations (AC1)--(AC3) only. We have obtained this presentation by 
employing a trick that was originally
 used (in a different context) by McCool and  Pietrowski \cite{McP}. 
Further examples of stably AC-equivalent sets that are not 
 known to be AC-equivalent are provided by our Theorem 1.4 below, but they 
are based on an altogether  different idea. 

 Also, we were able to show that every presentation in the series (4) above
 is AC-equivalent to a presentation in the series (5): 
 \smallskip 

 \noindent {\bf Proposition 1.3.} The presentation 
 $\langle x,y| ~x^n = y^{n+1}, ~xyx=yxy\rangle, 
  ~n \ge 2,$ is AC-equivalent to  
 $\langle x,y| ~x^{-1}y^n x = y^{n+1}, ~x=y^{-1}x^{-1}yxy\rangle$.
 \smallskip

  All these examples are discussed in Section 2. 
 In Section 3, motivated by the series (5) of balanced presentations of the trivial 
group,  
we give a similar generic series in the free group of rank 3. It is based  on 
a rather general and simple method, namely, on  obtaining balanced presentations  
of the trivial group from 
 ``exotic" knot diagrams of the unknot. 
This method can therefore produce numerous potential counterexamples to the
 Andrews-Curtis conjecture  in 
  a free group of any rank $\ge 2$. (Note however that these are {\it not} 
counterexamples to the {\it weak} Andrews-Curtis conjecture since it is known 
that Reidemeister moves  applied to a knot diagram give stably AC-equivalent 
Wirtinger presentations of the knot group.)  
In particular, we have: 
\smallskip 

\noindent {\bf Theorem 1.4.} Any presentation of the following form determines
the trivial group: $\langle  x, y, z ~| ~x = z \cdot [y^{-1}, x^{-1}, z], 
 ~y = x \cdot [y^{-1}, x^{-1}, z^{-1}] \cdot [z^{-1}, x], ~w \rangle$, ~where $w$  
is a word in 
$x$, $y$   and  $z$ whose exponent sum   on  $x$, $y$ and $z$  equals $\pm 1$.  
\smallskip 

 We note that the normal closure of the  first two relators contains the 
 commutator subgroup of the ambient free group $F_3$ since the corresponding
 two-relator group is the fundamental group of the unknot, i.e., is infinite
 cyclic. Hence, by  a sequence of AC-transformations, one can reduce 
 $w$ to a primitive element of $F_3$, and  therefore reduce the 3-generator
 presentation to a 2-generator one. For example, if in the  presentation in
Theorem 1.4, we 
 take $w=x^{-1}y z$, then, after simplifications, we get the following 
 2-generator presentation with the total length of relators equal to 25: 
 
 \noindent  $\langle  x, y ~| ~x^{-1} y^{-1} x y^{-1} x^{-1} y x y^{-2} x y x^{-1} y, 
~y^{-1} x^{-1} y^2 x^{-1} y^{-1} x y x  y^{-2} x \rangle$. 

Our computer program was able
to reduce  the total length of relators to 14: 
$$\langle  x, y ~| ~x y x^{-2} y^{-1} x y^{-1}, 
~x^{-1} y^{-1} x y^2 x y^{-1} \rangle,$$  
  but this latter presentation eludes all attempts on decreasing the  length.

 We also note that the method of constructing presentations of the trivial
 group based on ``exotic" knot diagrams of the unknot, is very flexible: 
 by changing various parameters of the method, one can obtain essentially 
 different presentations from the same knot diagram -- see the remark 
 after the proof of Theorem 1.4. Moreover, there is a lot of freedom 
 in choosing the word $w$, and  therefore in constructing potential counterexamples
of rank $n \ge 2$ based on generic presentations of rank $n + 1$. 
\medskip 

    In Section 4, we give yet another procedure for 
constructing  balanced presentations  
of the trivial group 
by using a composition method to   generate new balanced presentations
from old ones. We realize that this procedure may be known, 
 but we were unable to find it in the literature.

Briefly, the construction is as   follows. Let 
$$
P = \langle x_1, \ldots, x_n \mid r_1(X), \ldots , r_n(X) \rangle = \langle X
\mid R \rangle ,
$$
$$
Q = \langle x_1, \ldots, x_n \mid s_1(X), \ldots , s_n(X) \rangle = \langle X
\mid S \rangle
$$
(where $r_i(X), ~s_i(X)$ are words in generators $x_1, \ldots,x_n$) be two 
balanced presentations of the trivial group. We assume here that relators form
an ordered set, so that reordering   the relators results in a different
presentation. Define a new presentation $P\circ Q$ as follows: 
$$
P \circ Q = \langle x_1, \ldots, x_n \mid r_1(S), \ldots , r_n(S) \rangle = \langle X
\mid R(S) \rangle,
$$
where $r_i(S)$ is the word obtained from $r(X)$ upon replacing $x_j$ by $s_j$
for all $j$.

 Then, if $P = \langle X \mid R \rangle$ and   $Q = \langle X
\mid S \rangle$ were balanced presentations of the trivial group, the 
new presentation 
$P \circ Q$ is that of the trivial group as well (see Lemma 4.1).

 It is easy to
see that the set of all balanced presentations on generators $x_1, \ldots,
x_n$ forms a semigroup ${\cal P}_n$ with respect to the composition $\circ$.
 Clearly, ${\cal P}_n$ is embeddable into the semigroup of all
endomorphisms of $F_n$. Note that every finite subset $K \subset {\cal
P}_n$ generates a series of new presentations
-- those in the semigroup $\langle K \rangle$ generated by $K$ in ${\cal P}_n.$
This allows one to produce easily new examples of balanced presentations of
the trivial group. The following result shows that the Andrews-Curtis conjecture holds for
every presentation from ${\cal P}_n$ if and only if it holds for every
presentation from some generating set of ${\cal P}_n.$
\smallskip

\noindent {\bf Proposition 1.5.} (cf. \cite{MyasnikovJr}) Let $K$ be a subset 
of ${\cal P}_n$. Then every presentation in 
  the semigroup $\langle K \rangle$ generated by $K$ in ${\cal P}_n$
satisfies the Andrews-Curtis conjecture if and only if every presentation 
in $K$ does. 

\smallskip

    In Section 5, we consider the Andrews-Curtis equivalence in   
 metabelian groups of the form $F/[R,R]$, where $R$  is a 
normal subgroup of a free group $F$. (The definition of the 
AC-equivalence generalizes obviously to arbitrary groups.) 
We shall denote elements of  a  free  group  and  their  
 natural  images in a  group $F/[R,R]$ by the same letters when
there is  no ambiguity. 

Myasnikov \cite{Myasnikov} settled  the analog of the Andrews-Curtis
conjecture for  those  groups in the positive. More precisely, 
 he showed  (in the case where $F/R$ is a free abelian group of 
rank $n \ge 2$) that  if, in a   metabelian group $M=F/[R,R]$ ,
there are two sets of  cardinality $n$ whose normal closures equal
$M$, then those  sets are  AC-equivalent.  (In fact, his result is 
somewhat more general and  remains valid, in particular, for free 
 solvable groups of arbitrary derived length.)

 In sharp contrast, we prove here the following 
\smallskip 

\noindent {\bf Theorem 1.6.} Let $F=F_n$ be the free group of rank $n
\ge 3$ generated by  $X = {\{}x_1, ..., x_n {\}}$,  and
let  $R$ be the normal closure of $Y = {\{}[x_1, x_2], x_3 {\}}$. 
Then, in the  metabelian group $M=F/[R,R]$, there are infinitely
many pairwise AC-inequivalent, but stably  AC-equivalent 2-element
sets each of which has the same  normal closure $R$ in the group
$M$. A particular pair of elements with  the same  normal closure  as $Y$ but   AC-inequivalent to $Y$, is  ${\{}[x_1, x_2][x_1, x_2, x_1]^{-2x_2^{-1}} \cdot x_3^{4x_2^{-1}},
  ~[x_1, x_2, x_1, x_1]^{-x_2^{-1}} \cdot x_3[x_3, x_1]^{2x_2^{-1}}{\}}$. 
\smallskip 

 We note that  Lustig \cite{Lustig} and Metzler \cite{Metzler1}  
 constructed (finitely many) 
pairwise AC-inequivalent   sets  of the same cardinality ($\ge 3$) in  a free 
group,  each of which has the same  normal closure.
 It is not known however if an {\it infinite} family of sets like
that  exists in a free group. Furthermore,  the sets constructed 
in \cite{Lustig} and \cite{Metzler1} are also stably  
AC-inequivalent since the corresponding 2-complexes are
simple-homotopy  inequivalent. It seems to be another very difficult 
question  as to  whether or not there are stably
AC-equivalent,  but AC-inequivalent sets of the same cardinality in
a free group. 

 It is therefore amazing that metabelian groups  not only admit
an easier solution of the problem analogous to the Andrews-Curtis
conjecture, but also admit easier counterexamples to its 
generalization. 

 Our proof of Theorem 1.6 is based on the fact that there are 
non-tame invertible 
$2\times 2$ matrices over Laurent polynomial  rings
in 2 or more variables. (An invertible square matrix is called {\it tame}
if it is a product of  elementary and  diagonal  matrices). This fact is
due to Bachmuth and Mochizuki \cite{BM}.   A particular matrix like that was recently found by  Evans \cite{Evans}. 

On the other hand, 
 every invertible square matrix over a Laurent  polynomial  ring is {\it stably
tame}, which means  one can extend this matrix  by placing 1 on the
diagonal and 0 elsewhere, such that the extended matrix is tame. This
follows from Suslin's stability  theorem \cite{Suslin}.

    In  view  of the aforementioned   results   of   Suslin   and 
Bachmuth-Mochizuki,  the  following  open  problem is of particular 
interest (see \cite{Ba} for a survey): 
{\it is it true that every  matrix  from} $GL_2({\bf Z}[t^{\pm 1}])$ {\it is 
tame?}  Here  ${\bf Z}[t^{\pm 1}]$  denotes the integral group ring  of  the 
infinite  cyclic  group, i.e., the ring of Laurent polynomials in one
variable.   We  show  that this notorious problem in $K$-theory 
is also related to the AC-equivalence: 
\smallskip 

 \noindent {\bf Proposition 1.7.} 
Every matrix from $GL_2({\bf Z}[t^{\pm 1}])$  is tame if and only if 
 every  set $\widetilde Y$ such that  $<\widetilde Y>^M =<Y>^M  = R$,   is
AC-equivalent to $Y={\{}x_1, x_2{\}}$ in the group $M=F_3/[R,R]$.\\

\noindent {\bf 2. AC-equivalent presentations} 
\bigskip

 We start by listing two  extra  transformations that  facilitate computations 
  toward establishing the AC-equivalence of various sets. Either of them 
 is a composition of several transformations (AC1)--(AC3). 
  The first transformation 
 comes from \cite{Olga}, where it is called ``Basic substitution 
principle". 
\smallskip 

\noindent {\bf Making a substitution.} If in a set $Y = {\{}y_1, ..., y_m
{\}}$, one replaces some element $y_i$ by an element $\tilde y_i$ which
is  congruent to $y_i$ modulo the normal closure of $Y = {\{}y_1, ...,
y_{i-1}, y_{i+1}, ..., y_m {\}}$, then the resulting set  $\widetilde  Y$ 
will be AC-equivalent to $Y$. 
\smallskip 

 The second transformation is a composition of several  (AC1)--(AC3) only  if 
 the normal closure of a given set $Y$ is the whole group  $F$. 

\smallskip 

\noindent {\bf Applying an automorphism.} If $\langle Y\rangle^F  = F$, then
 applying any automorphism of 
the ambient free group $F$ to {\it every element} of  the set $Y$ 
gives a set which is AC-equivalent to $Y$. 
\smallskip 

 This last statement has the following useful corollary:
\smallskip 

\noindent {\bf Corollary 2.1.} If $\langle x,y| ~r_1, r_2\rangle$ is a 
presentation of 
the trivial group  and $r_1$ is a primitive element of the ambient free group
$F_2$, then this presentation is AC-equivalent to the  standard one 
$\langle x,y| ~x, y\rangle$. 
\smallskip 

 Now we are ready for a proof of Proposition 1.1. We shall just give 
a relevant chain of AC-transformations in each case (skipping obvious
steps sometimes), emphasizing the 
most subtle steps  in a chain by placing a 
(!)  after them.  (This is influenced by   chess notation.) 
\smallskip 

\noindent {\bf Proof of Proposition 1.1.} 

\noindent {\bf (a)} $\langle x,y| ~x^2 = y^3, xyx=yxy\rangle
\longrightarrow \langle x,y| ~x^{-1}y^3x^{-1}, yxy^{-1}x^{-1}y^{-1}x\rangle 
\longrightarrow \\
\langle x,y| ~x^{-1}y^3x^{-1}, yxy^{-1}x^{-1}y^2x^{-1}\rangle 
\longrightarrow \langle x,y| ~yx^{-1}y^3x^{-1}y^{-1}, 
yxy^{-1}x^{-1}y^2x^{-1}\rangle \\
  \longrightarrow 
  \langle x,y| ~yx^{-1}y^2x^{-1}y^2x^{-1}, ~yxy^{-1}x^{-1}y^2x^{-1}\rangle$.

 Now the first relator is  a primitive element (say, the automorphism 
 $x \to xy^2; ~y \to y$ ~takes it to $y^{-1}x^{-3}$), therefore we are
 done. 
\smallskip 

\noindent {\bf (b)} $\langle x,y| ~x^{-1}yx =y^2, x=yx^2yx^{-2}\rangle 
\longrightarrow \langle x,y| ~x^{-1}yx =y^2, x^{-1}yx^2yx^{-2}\rangle 
\longrightarrow $ (substitution of $y^2$ for $x^{-1}y x$) 
$~\langle x,y| ~x^{-1}yx =y^2, y^2xyx^{-2}\rangle 
\longrightarrow \\
\langle x,y| ~x^{-1}yx =y^2, x^{-1}y^2xyx^{-1}\rangle 
\longrightarrow $ (substitution of $y^4$ for $x^{-1}y^2 x$) \\
$~\langle x,y| ~x^{-1}yx =y^2, y^5x^{-1}\rangle ~(!)$. 

 Now the  second relator is  a primitive element, therefore we are done. 

\smallskip 

\noindent {\bf (c)} $\langle x,y| ~x^{-1}y^2 x = y^3, x^2 = yxy^{-1}\rangle
\longrightarrow \langle x,y| ~y^{-2}x^{-1}y^2 xy^{-1}, x^{-2}yxy^{-1}\rangle
\longrightarrow 
\langle x,y| ~y^{-2}x^{-1}y^2 xy^{-1}, yx^{-1}y^{-1}x^2\rangle 
\longrightarrow \langle x,y| ~y^{-2}x^{-1}yx^2, yx^{-1}y^{-1}x^2\rangle 
\longrightarrow \\
\langle x,y| ~xy^{-2}x^{-1}yx, xyx^{-1}y^{-1}x\rangle 
\longrightarrow \langle x,y| ~xyxy^{-2}x^{-1}yxy^{-1}x^{-1}, 
xyx^{-1}y^{-1}x\rangle \\
  \longrightarrow 
\langle x,y| ~xyxy^{-2}, xyx^{-1}y^{-1}x\rangle 
\longrightarrow \langle x,y| ~y^{-1}xyxy^{-1}, yx^{-1}y^{-1}x^2\rangle 
\longrightarrow \\
\langle x,y| ~y^{-1}x^3, yx^{-1}y^{-1}x^2\rangle$. 

\smallskip 

\noindent {\bf (d)} $\langle x,y| ~x^{-1}y^2 x = y^3, ~x^2 = yxy\rangle
\longrightarrow \langle x,y| ~x^{-1}y^2 x = y^3, xy^{-1}=x^{-1}yx\rangle 
\longrightarrow$ (substitution of $(xy^{-1})^2$ for $x^{-1}y^2 x$) 
$~\langle x,y| ~xy^{-1}xy^{-1} = y^3, xy^{-1}=x^{-1}yx\rangle (!) \\
\longrightarrow$  (applying the automorphism $x \to xy; ~y \to y$) 
$\langle x,y| ~x^2 = y^3, x=y^{-1}x^{-1}yxy\rangle 
\longrightarrow 
\langle x,y| ~x^2 = y^3, xyx=yxy\rangle$. 

  Now the result follows from part (a). $\Box$ 

\medskip

\noindent {\bf Proof of Proposition 1.2.} 
  By  Proposition 1.3, $\langle x,y| ~x^3 = y^4, ~xyx=yxy\rangle$ is 
AC-equivalent to $\langle x,y| ~x^{-1}y^3x=y^4, ~xyx=yxy\rangle$, so we can start from 
there: 

\noindent $\langle x,y| ~x^{-1}y^3x=y^4, ~xyx=yxy\rangle 
\longrightarrow$ ($x \to x^{-1}; ~y \to y^{-1}$) 

\noindent $\langle x,y| ~x^{-1}y^{-1}x^{-1}=y^{-1}x^{-1}y^{-1}, 
~xy^{-3}x^{-1}=y^{-4}\rangle \longrightarrow$ 

\noindent $\langle x,y| ~xyx=yxy, ~xy^{-3}x^{-1}=y^{-4}\rangle 
\longrightarrow$ ($x \to xy^{-1}; ~y \to y$) 

\noindent $\langle x,y| ~x^2y^{-1} =
yx, ~xy^{-3}x^{-1}=y^{-4}\rangle \longrightarrow \langle x,y| ~y^{-1}x^2=xy, 
~xy^{-3}=y^{-4}x\rangle \longrightarrow$ 
 $\langle x,y| ~y^{-1}x^2=xy, 
~xy^{-3}=y^{-4}xyy^{-1}\rangle \longrightarrow$ (substitution) \\ 
$\langle x,y| ~y^{-1}x^2=xy, 
~xy^{-3}=y^{-5}x^2y^{-1}\rangle (!) \longrightarrow \langle x,y| ~y^{-1}x^2=xy, 
~y^{-5}x^2=xyy^{-3}\rangle \longrightarrow$ 
(substitution)   $\langle x,y| ~y^{-1}x^2=xy, 
~y^{-4}x^2=x^2y^{-3}\rangle \longrightarrow \\
\langle x,y| ~x^2=yxy, ~y^{-4}x^2=x^2y^{-3}\rangle 
\longrightarrow $

\noindent $\langle x,y,z| ~x^2=yxy, ~y^{-4}x^2=x^2y^{-3}, ~z=y^{-3}\rangle 
\longrightarrow $

\noindent $\langle x,y,z| ~x^2=yxy, ~y^{-1}=x^2y^{-3}x^{-2}y^3, ~z=y^{-3}\rangle 
\longrightarrow $

\noindent $\langle x,y,z| ~x^2=yxy, ~y^{-1}=x^2 z x^{-2} z^{-1}, 
~z=(x^2 z x^{-2} z^{-1})^3\rangle \longrightarrow $

\noindent $\langle x,y,z| ~x^2=z x^2 z^{-1} x^{-1} z x^2 z^{-1} x^{-2}, 
~y=z x^2 z^{-1} x^{-2}, ~z=[x^2, z]^3 \rangle \longrightarrow $

\noindent $\langle x,z| ~x^4= z x^2 z^{-1} x^{-1} z x^2 z^{-1}, 
~z=[x^2, z]^3\rangle  \longrightarrow$ ($z \to y$)  

\noindent $\langle x,y| ~x^4= y x^2 y^{-1} x^{-1} y x^2 y^{-1},  
~y=[x^2, y]^3\rangle$. $\Box$ 
\medskip

\noindent {\bf Proof of Proposition 1.3.} 
  $\langle x,y| ~x^n = y^{n+1}, ~xyx=yxy\rangle 
\longrightarrow \langle x,y| ~x^n = y^{n+1}, ~x=(y^{-1}x^{-1})y(xy)\rangle 
\longrightarrow $ (substitution) 
$~\langle x,y| ~y^{-1}x^{-1}y^n xy = y^{n+1}, ~x=(y^{-1}x^{-1})y(xy)\rangle (!)
~\longrightarrow \langle x,y| ~x^{-1}y^n x = y^{n+1}, 
~x=y^{-1}x^{-1}yxy\rangle$.  $\Box$ \\

\noindent {\bf 3. Potential counterexamples from knot diagrams}
\bigskip

  In this section,  we give a rather general and simple method for constructing  
  potential counterexamples  to the original Andrews-Curtis 
 conjecture based on  ``exotic" knot diagrams 
of the unknot. Note however  that it cannot produce counterexamples to the weak 
 Andrews-Curtis conjecture 
 since it is known that Reidemeister moves  applied to a knot diagram give 
stably AC-equivalent Wirtinger presentations of the knot group.

    The method  can be  summarized as follows.
  Take a piece of rope, glue the ends, crumple (or even tie) 
  the rope, and through it on the table. What you have now on the table,
  is an ``exotic" diagram of the unknot. (We note that you can get a 
  diagram of the unknot by just drawing a closed curve on a list of paper, 
  making each self-crossing an undercrossing until you eventually close
  the curve, but this method  cannot possibly give you a desired counterexample!) 
  
 Now comes the mathematical part of the method. If you are lucky, the number 
of crossings in your knot diagram cannot be reduced by a single 
Reidemeister move (see e.g. \cite[p.9]{Kauf}). If it can, you have to  
``manually" adjust your diagram 
until you get one with that property (otherwise, you will have to do   
extra work in the next stage).   Now    write down the Wirtinger presentation 
of the fundamental group of the 
unknot (which is the infinite cyclic group) based on your knot diagram. This is 
a balanced presentation, with the number of generators equal to the 
number of crossings in  the knot diagram. (A particularly simple explanation 
of how to get the Wirtinger presentation from a knot diagram is given on 
p.268 of the book \cite{Kauf}).  
If the ambient free group is generated by 
$x_1, ..., x_n$, then the normal closure of the relators is the same as that 
of $x_1x_2^{-1}, ..., x_1x_n^{-1}$. All  Wirtinger 
relators are of the form  $x_i=x_j x_k x_j^{-1}$, or $x_i=x_j^{-1} x_k x_j$.
Furthermore,  any one of the relators is 
redundant, i.e., it follows from the other $(n-1)$. 

 Based on these well-known facts, we see that if we discard one relator
 from our Wirtinger presentation, and then add a relator   
which would trivialize the group, we shall get a balanced presentation 
of the trivial group.  We are now ready for 
\smallskip 

\noindent {\bf Proof of Theorem 1.4.} Consider the following diagram 
of the unknot (Figure 1).

\begin{figure}[htbp]
\centerline{\epsfig{file=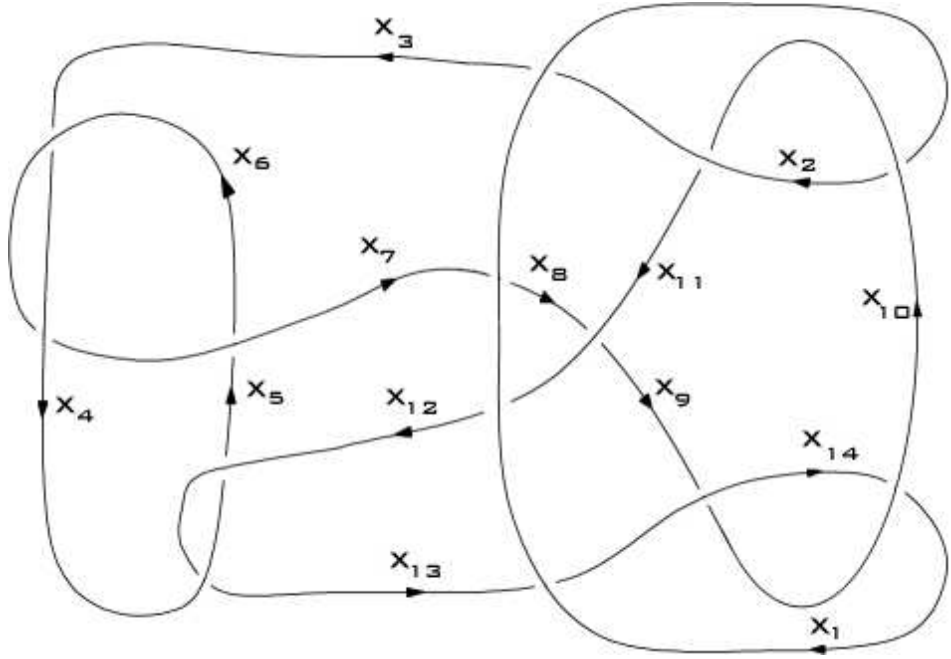}}
\caption{A diagram of the unknot.}
\label{fig1}
\end{figure}
 
 The corresponding Wirtinger presentation is 
 $\langle x_1, x_2, x_3, x_4, x_5, x_6, x_7, x_8,\\    
x_9, x_{10}, x_{11}, 
 x_{12}, x_{13}, x_{14} ~| ~x_1 = x_{10}x_{14}x_{10}^{-1}, 
 ~x_2=x_{10}^{-1}x_1x_{10}, ~x_3=x_1^{-1}x_2x_1,\\  
x_4=x_6^{-1}x_3 x_6, 
 ~x_5=x_{12}x_4x_{12}^{-1}, ~x_6=x_7^{-1}x_5x_7, ~x_7=x_4^{-1}x_6 x_4, 
~x_8=x_1x_7x_1^{-1}, ~x_9=x_{11}^{-1}x_8x_{11}, ~x_{10}=x_{14}x_9x_{14}^{-1},
~x_{11}=x_2^{-1}x_{10}x_2, ~x_{12}=x_1^{-1}x_{11}x_1, \\
 x_{13}=x_5x_{12}x_5^{-1},
 ~x_{14}=x_1x_{13}x_1^{-1} \rangle$.
\smallskip 

  Now we discard the relator $x_{12}=x_1^{-1}x_{11}x_1$, ~and then eliminate the 
generators one by one, starting with $x_{14}=x_2$ (from the first two relators).
That means, we apply relevant AC-transformations (cf. the ``basic
 substitution principle" in the previous section) to the first two relators 
 to get $x_{14}=x_2$; then replace $x_{14}$ by $x_2$ 
 ~in all relators that involve $x_{14}$ (the ``basic
 substitution principle" again), and then get rid 
of $x_{14}$ by applying the transformation (AC5).

 In a similar fashion, we successively   eliminate 
$~x_{11}=x_9, ~x_9=x_8$,   $~x_{13}, ~x_8, ~x_4, 
x_3, ~x_6, ~x_{10}, ~x_2$, 
 and $x_1$.  In the end, we get the following 
presentation: $\langle  x_5, x_7,  x_{12}~| ~x_5=x_{12}x_7^{-1}x_5^{-1}x_7x_5x_{12}
x_5^{-1}x_7^{-1}x_5x_7x_{12}^{-1}, ~x_7=x_5^{-1}x_7x_5x_{12}^{-1}x_5^{-1}
x_7^{-1}x_5x_7x_5x_{12}x_5^{-1}x_7^{-1}x_5 \rangle$.  

 Denote $x=x_5, ~y=x_7, ~z=x_{12}$. To make this group trivial, we have to 
add a relator   which would, together with the ones that we have, generate the
group $F_3=\langle  x, y, z\rangle$ modulo its commutator subgroup. 
 Since the abelianization of our relators is ${\{}x z^{-1}, x^{-1}y{\}}$, 
 adding any relator whose 
 exponent sum   on  $x$, $y$ and $z$  equals $\pm 1$, will trivialize the group. 
 This completes the proof. $\Box$ 
\smallskip 

\noindent {\bf  Remark.}  The presentation of the trivial group obtained by 
this method   {\it  very much} depends on the following: 

\noindent (1) the  relator   discarded in the very beginning; 

\noindent (2)  the order  in which generators are eliminated.  

By varying these choices, one can obtain several  essentially different 
balanced presentations of the trivial group from the same 
knot diagram.\\

 \noindent {\bf 4. The semigroup of balanced presentations of the trivial group}
\bigskip

 We start with 
 \smallskip

\noindent {\bf Lemma 4.1.} Let $P = \langle X \mid R \rangle,$ $Q = \langle X
\mid S \rangle$ be two balanced presentations of the trivial group. Then
$P \circ Q = \langle X \mid R(S) \rangle $ also presents the trivial group.
\smallskip

\noindent {\bf Proof.}  Since the group defined by $P$ is trivial, every
element $x_i \in X$ can be expressed  as a product of
conjugates of relators from $R,$ say
$$
x_i = \prod_j r_{ij}(X)^{w_{ij}(X)}.
$$
Upon replacing $x_i$ by $s_i(X)$ 
(i.e., upon applying an endomorphism of $F(X)$ to both sides), we get
$$
s_i = \prod_j r_{ij}(S)^{w_{ij}(S)}.
$$
 Therefore, $S \subset \langle R(S) \rangle ^F$, hence $\langle R(S) \rangle
^F \supseteq \langle S \rangle ^F = F.$ This shows that $P \circ Q$  defines
the trivial group. $\Box$
\smallskip

\noindent {\bf Proof of Proposition 1.5.} It suffices to show that if a balanced
presentation
$P = \langle X \mid R \rangle$ is AC-equivalent to the standard one, then for
any balanced presentation $Q = \langle X \mid S \rangle$ the composition $P
\circ Q = \langle X \mid R(S) \rangle $ is AC-equivalent to the $Q.$

Let $t_1,...,t_n$ be a sequence of transformations (AC1)-(AC3) that takes 
$P$ to the standard presentation of the trivial group. Denote by $u_1,...,u_n$
a chain of Andrews-Curtis transformations obtained from $t_1,...,t_n$ as
follows:
\smallskip

\noindent  -- if $t_i$ is a transformation of the type (AC1) or (AC2), then
$u_i=t_i$;\\ 
\noindent  --  if $t_i$ is of the type (AC3) (say, $t_i$ replaces $r_j$ by
$w(X)r_j w(X)^{-1}$), then $u_i$ replaces $r_j$ by $w(S)r_j w(S)^{-1}$.
\smallskip

 Obviously, the sequence $s_1, \ldots, s_n$ transforms $P \circ Q$ into $Q.$
$\Box$

\smallskip

\noindent {\bf   Example.} Let $P = \langle x,y \mid y^{-1} x y = x^2, 
x^{-1} y x = y^2 \rangle.$ 
It is easy to see that $P$ is AC-equivalent to the
standard presentation $\langle X \mid X \rangle.$ If $Q = \langle x,y \mid r,s
\rangle$ is an arbitrary presentation of the trivial group, then, by 
Proposition 1.5, the presentation
$$
P\circ Q = \langle x,y \mid r^{-1} s r= s^2, s^{-1} r s = r^2 \rangle
$$
is AC-equivalent to $Q.$ \\

\noindent {\bf 5. The Andrews-Curtis equivalence in metabelian
groups} 
\bigskip

    We  recall that given a group  $G$  presented in the  form   $G  = 
F/R$,  one can turn the abelian group  $R/R'$  into a left  ${\bf Z}G$-module  upon 
setting  $g(rR') = g^{-1} r g R',  ~g \in G, ~r \in  R$,  and extending this 
 action 
${\bf Z}$-linearly to the whole ring ${\bf Z}G$. 

    One more  definition  that  we  need  seems  to  be  well-known:  a 
presentation  $G  = F/R$ is  called  aspherical  if  the  corresponding 
relation module  $R/R'$   is  a free ${\bf Z}G$-module.  Now we have 
\smallskip 

 \noindent {\bf Proposition 5.1.} Let a group $G$  be given by an aspherical 
 presentation  $G  = F/R$,   and let the corresponding relation module $R/R'$ 
be  generated  by  $k$  elements.  Then: any two generating systems of the 
relation module $R/R'$ are  (stably) AC-equivalent if and only if every matrix
from the  group   $GL_k({\bf Z}G)$   is (stably)  tame,  and  all  diagonal 
matrices   in   the   corresponding decompositions have trivial units
(i.e., elements of the form  $\pm g,  ~g \in G$) on the diagonal. 
\smallskip 

    This suggests two ways of constructing counterexamples.   The  first 
way is to pick up an aspherical presentation  $G  = F/R$ of  a  group  
$G$  such that the group ring ${\bf Z}G$ has non-trivial units;  this  idea  has 
been essentially used in \cite{Lustig}. The second way, the one that we use
 here, is to choose  a  group $G$  such that there are non-tame invertible 
 matrices over the group ring ${\bf Z}G$.
\medskip

 \noindent {\bf Proof of Proposition 5.1.} Fix a generating system  
 $Y = {\{}r_1R',...,r_kR'{\}}$  of the 
relation module $R/R'$, and  let   $Y  =  {\{}s_1R',...,s_kR'{\}}$   be  another 
generating system.  Then we can write: 
\begin{equation}
s_i \equiv  \prod_{j=1}^{k} r_j^{u_{ij}}    ~(mod ~R'),  ~1 \le i \le  k,
\end{equation}                
where  $u_{ij} \in {\bf Z}G$.  Similarly, 
\begin{equation}
r_i \equiv  \prod_{j=1}^{k} s_j^{v_{ij}}    ~(mod ~R'),  ~1 \le i \le  k,
\end{equation}  
for some  $v_{ij} \in {\bf Z}G$.  Hence, if we  define two $k \times k$ matrices   $U = (u_{ij})$         
and  $V = (v_{ij})$  over the group ring  ${\bf Z}G$,  then (1) and (2) can be 
written in the matrix form as  $(s_1R',...,s_kR')^t  = 
U \cdot (r_1R',...,r_kR')^t$, 
and  $(r R',...,r R')^t  = V \cdot  (s R',...,s R')^t$,  where $^t$  means  transposition, 
i.e., we consider our $k$-tuple as a column.  The multiplication by elements
 of ${\bf Z}G$ is  interpreted  here as the module action described above. 
 
This yields  
$VU \cdot (r_1R',...,r_kR')^t  = (r_1R',...,r_kR')^t$.  Since  
$R/R'$  is  a  free ${\bf Z}G$-module and  ${\{}r_1R',...,r_kR'{\}}$  a free basis, this is equivalent to  
$VU = I$,  the identity matrix,  therefore 
\begin{equation}
                             U, V  \in GL_k({\bf Z}G).         
\end{equation}                
     We now observe that applying an elementary transformation (AC1)  to 
a $k$-tuple \\
$(r_1R',...,r_kR')$  amounts  to  multiplying  the  corresponding 
column by an elementary matrix on the left, and applying (AC2) or  (AC3) 
amounts to multiplying it by a diagonal matrix which  has  only  trivial 
units on  the  diagonal.   Conversely,  a multiplication  by  any  such  a 
diagonal or elementary matrix on  the  left  can  be  interpreted  as  a 
(sequence of) elementary transformation(s) (AC1)--(AC3).  Combining  this 
with (3) yields the result. $\Box$ 
\smallskip 

\noindent {\bf Proof of Theorem 1.6} follows almost immediately
from Proposition 5.1. The group $G$ in this situation is the 
free abelian group of rank 2 generated by the natural images of $x_1$ and $x_2$; the presentation $<x_1,x_2, x_3| 
~[x_1, x_2], x_3 >$ is obviously aspherical. The group ring ${\bf Z}G$ 
is the ring of Laurent polynomials in 2  variables. Over this 
ring, there are non-tame $2 \times 2$ invertible matrices -- see 
\cite{BM}. Moreover, the index of the subgroup $GE_2({\bf Z}G)$ 
generated by elementary and  diagonal  matrices, in the whole group
$GL_2({\bf Z}G)$, is infinite \cite{BM}. Hence, the first part of 
Theorem 1.4 follows now from Proposition 5.1. The particular pair of elements  
in the statement of Theorem 1.6 is modeled on the example of a matrix 
 in $GL_2({\bf Z}G)$, but not in $GE_2({\bf Z}G)$, given in \cite{Evans}.  
The matrix is: 

  $\left(\begin{array}{cc} 1-2(x_1-1)x_2^{-1} & 4x_2^{-1}\\ 
 -(x_1-1)^2x_2^{-1} & 1+2(x_1-1)x_2^{-1}  \end{array}\right)$\\    
 
 On the other hand, by a result of Suslin \cite{Suslin}, 
$GE_k({\bf Z}G)=GL_k({\bf Z}G)$ for any $k \ge 3$; in particular, every 
matrix from $GL_2({\bf Z}G)$ is 
stably tame, whence the  result about stable AC-equivalence. $\Box$ 
\medskip 

\noindent {\bf  Proof of  Proposition 1.7.} In this situation,
the group $G$ is infinite cyclic (generated by $x_3$), so that
the group ring ${\bf Z}G$ is the ring of one-variable Laurent polynomials.
 The result now follows from Proposition 5.1. $\Box$ \\

\noindent {\bf Acknowledgement} 
\medskip 

 We are grateful to   M. Lustig for useful  discussions.

\medskip
\noindent 
 Department of Mathematics, The City  College  of New York, New York, 
NY 10031 
\smallskip

\noindent {\it e-mail addresses\/}:  alex@rio.sci.ccny.cuny.edu, 
alexei@rio.sci.ccny.cuny.edu,\\  
\noindent shpil@groups.sci.ccny.cuny.edu 

\end{document}